\documentclass[12pt,a4paper]{article}
\RequirePackage{amsmath}%
\RequirePackage{amsthm}%
\usepackage{amsfonts}%
\usepackage{amssymb}%
\usepackage{graphicx}%

\newtheorem{theorem}{Theorem}
\newtheorem{lemma}{Lemma}
\newtheorem{corollary}{Corollary}

\newcommand{\codim}{{\rm codim}}
\newcommand{\Range}{{\rm Range}}
\newcommand{ \Null}{{\rm Null}}

\renewcommand{\Re}{\mathop{\mathrm{Re}}}
\renewcommand{\Im}{\mathop{\mathrm{Im}}}

\renewcommand{\i}{\mathrm{i}}

\newcommand{\bI}{{\bf I}}

\newcommand{\bM}{{\bf M}}
\newcommand{\bN}{{\bf N}}

\begin{document}

\title{The Riemann-Hilbert problem \\and the generalized Neumann kernel \\on
unbounded multiply connected regions\footnote{This paper has been published in:
The University Researcher (IBB University Journal), 20 (2009) 47--60.}}

\author{ Mohamed M.S. Nasser  }

\date{}
\maketitle

%\centerline{ {\tt Draft1.tex}}
\vskip-1.0cm %
\centerline{Department of Mathematics, Faculty of
Science, Ibb University,} %
\centerline{P.O.Box 70270, Ibb, Yemen}%
\centerline{E-mail: mms\_nasser@hotmail.com}

\vskip0.7cm

\begin{center}
\begin{quotation}
{\noindent {\bf Abstracts.\;\;}%
A Fredholm integral equation of the second kind with the generalized Neumann kernel associated with the Riemann-Hilbert problem on unbounded multiply connected regions will be derived and studied in this paper. The derived integral equation yields a uniquely solvable boundary integral equations for the modified Dirichlet problem on unbounded multiply connected regions.
}%
\end{quotation}
\end{center}
\begin{center}
\begin{quotation}
\end{quotation}
\end{center}
\begin{center}
\begin{quotation}
\end{quotation}
\end{center}
%

%-----------------------INTRO--------------------------------------
\section{Introduction}
\label{sect:intr}
%------------------------------------------------------------------

%
The interplay of Riemann-Hilbert problems (RH problems, for short)
and integral equations with the generalized Neumann kernel has
been investigated in~\cite{wegg} for simply connected regions,
in~\cite{wegm} for bounded multiply connected regions, and
in~\cite{nasr} for simply connected regions with piecewise smooth
boundaries.
Integral equations with the generalized Neumann kernel on bounded
multiply connected regions have been used in~\cite{nasc} to develop
a unified method for computing conformal mapping onto the
classical slit domains.
Based on the results of this paper, the method presented
in~\cite{nasc} can be extended to unbounded multiply connected
regions (see~\cite{nasm}).
In this paper, we shall extend the results of~\cite{wegg,wegm} to
unbounded multiply connected regions.
We derive a second kind Fredholm integral equation with the
generalized Neumann kernel for the RH problems on unbounded
multiply connected regions.
Then, based on a M\"obius transform, the properties of the derived
integral equation will be deduced from the properties of integral
equations with the generalized Neumann kernel on bounded multiply
connected regions.
The derived integral equation yields a boundary integral equation
for the modified Dirichlet problem which is a special case of the
RH problem.
The remainder of this paper is organized as follows: we present
some auxiliary material in Section~\ref{sc:aux}.
Section~\ref{sc:gnk} presents an integral equation with the
generalized Neumann kernel for the RH problem on unbounded
multiply connected regions.
The solvability of the derived integral
equation will be studied
in Section~\ref{sc:slv}.
Section~\ref{sc:mdir} presents a uniquely
solvable integral equation for the modified Dirichlet problem.
Finally, short conclusions will be given in Section~\ref{sc:con}.

%--------------------------------------------------------------
\section{Auxiliary material}
\label{sc:aux}
%--------------------------------------------------------------

\begin{figure}[!h]
\centerline{\scalebox{0.6}[0.3]{
\includegraphics{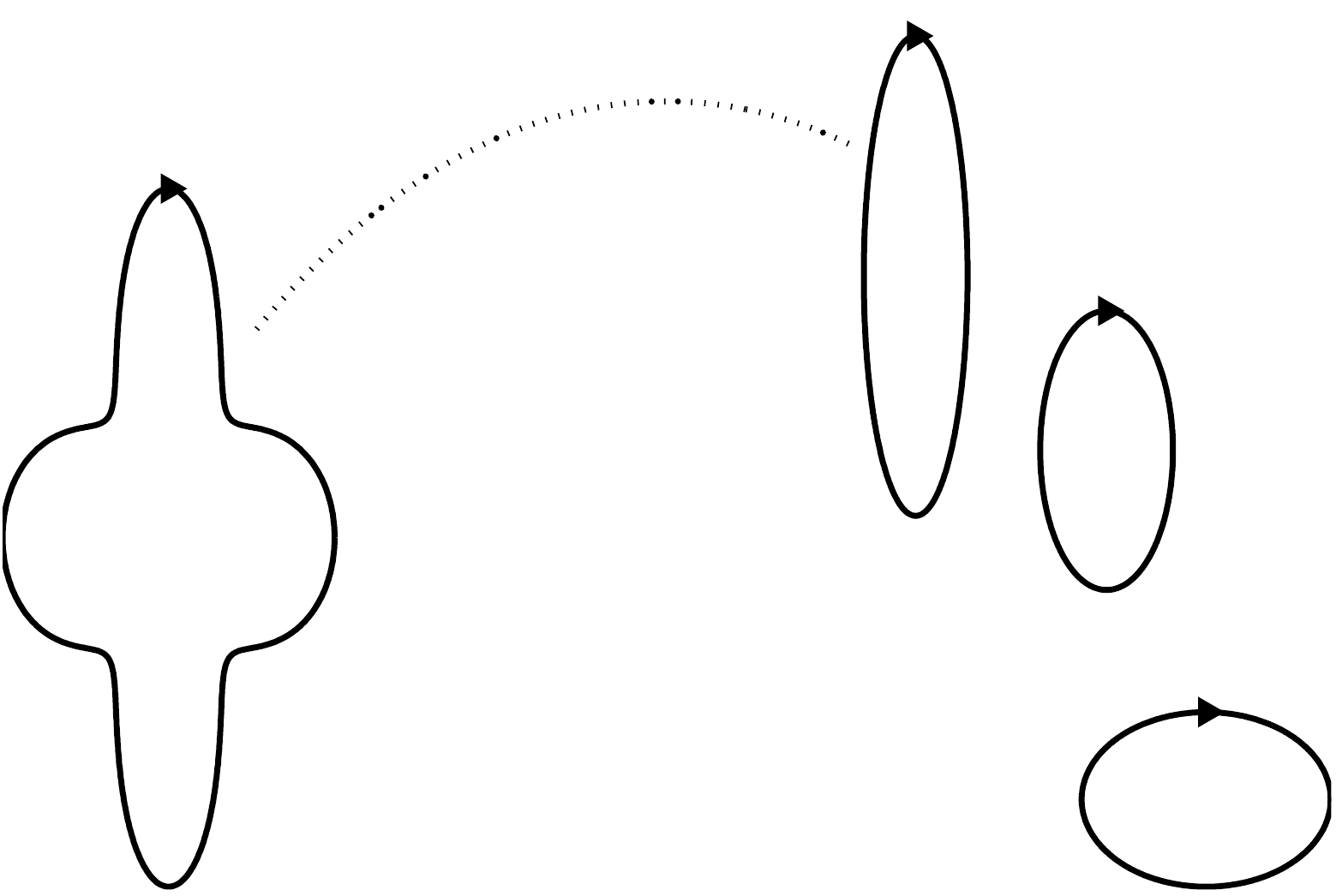}}}
 \vskip-3.7cm \noindent\hspace{8.95cm} $\Gamma_3$
 \vskip-.0cm \noindent\hspace{3.6cm} $\Gamma_m$
 \vskip-.1cm \noindent\hspace{10.25cm} $\Gamma_2$
 \vskip-0.3cm \noindent\hspace{8.9cm} $G_3$
 \vskip+.1cm \noindent\hspace{6.9cm}$G$ \hspace{3.0cm}$G_2$
 \vskip-0.2cm \noindent\hspace{3.5cm} $G_m$
 \vskip-.2cm  \hspace{10.95cm} $\Gamma_1$
 \vskip+0.1cm \noindent\hspace{10.95cm} $G_1$
 \caption{\rm An unbounded multiply connected region $G$ of
connectivity $m$.}
 \label{f:unbd}
\end{figure}

Suppose that $G$  is an unbounded multiply connected region of
connectivity $m$ bounded by $\Gamma:=\partial G=
\Gamma_1\cup\Gamma_2\cup\cdots\cup\Gamma_m$ where the curves
$\Gamma_k$, $k=1,2, \ldots, m$, are simple non-intersecting smooth
clockwise oriented closed curves (see Figure~\ref{f:unbd}).
We assume that $\infty\in G$ and $0\in G$.
The complement $G^-:=(\mathbb{C}\cup\{\infty\})\backslash (G\cup
\Gamma)$ of $G$ consists of $m$ bounded simply connected regions
$G_k$, $k=1,2,\ldots,m$.
The curve $\Gamma_k$ is parametrized by a $2\pi$-periodic twice
continuously differentiable complex function $\eta_k(s)$ which
transverses $\Gamma_k$ in the clockwise orientation with $\dot
\eta_k(s)=d\eta_k/ds \ne 0$.
The total parameter domain $J$ is the disjoint union of $m$
intervals $J_k=[0,2\pi]$, $k=1,2,\ldots,m$.
We define a parametrization of the whole boundary as the complex
function $\eta$ defined on $J$ by
\begin{equation}\label{e:eta}
\eta(s)=\left\{ \begin{array}{l@{\hspace{0.5cm}}l}
\eta_1(s),&s\in J_1=[0,2\pi],\\
\hspace{0.3cm}\vdots\\
\eta_m(s),&s\in J_m=[0,2\pi].
\end{array}
\right .
\end{equation}

%
%%ggggggggggggggggggggggggggggggggggggggggggg
Let $H$ be the space of all real H\"older continuous
$2\pi-$periodic functions $\phi(s)$ of the parameter $s$ on $J_k$
for $k=1,2, \ldots, m$, i.e.,
\[
\phi(s) = \left\{
\begin{array}{l@{\hspace{0.5cm}}l}
 \phi_1(s),     & s\in J_1, \\
  \vdots,       & \\
 \phi_m(s),     & s\in J_m, \\
\end{array}%
\right.
\]
where $\phi_1,\ldots,\phi_m$ are real H\"older continuous
$2\pi-$periodic functions.
In view of the smoothness of $\eta$, a real H\"older continuous
function $\hat\phi$ on $\Gamma$ can be interpreted via $\phi(t) :=
\hat\phi(\eta(t))$ as a function $\phi\in H$; and vice versa.
If $\phi\in H$ is a piecewise constant real-valued function, i.e.,
\[
\phi(s) = \left\{
\begin{array}{l@{\hspace{0.5cm}}l}
 \alpha_1,     & s\in J_1, \\
  \vdots,       & \\
 \alpha_m,     & s\in J_m, \\
\end{array}%
\right.
\]
with real constants $\alpha_1, \ldots, \alpha_m$, then $\phi$ will
be denoted by
\[
\phi = (\alpha_1, \ldots, \alpha_m).
\]
Let $A$ be a continuously differentiable complex function on
$\Gamma$ with $A\ne 0$.
We assume that $A$ is given in the parametric form $A(t)$ such
that $t\mapsto A(t)$ is continuously differentiable for all $t\in J$.
With $\gamma,\mu\in H$, we define the function
\begin{equation}
\label{e:Phi} \Phi(z):=\frac{1}{2\pi \i} 
\int_\Gamma \frac{\gamma+\i\mu}{A} \frac{d\eta}{\eta-z}, \quad z\notin \Gamma.
\end{equation}
Then $\Phi$ is an analytic function in $G$ as well as in $G^-$
with $\Phi(\infty)=0$. The boundary values $\Phi^+$ from inside
$G$ (the left of $\Gamma$) and $\Phi^-$ from outside $G$ (the
right of $\Gamma$) are H\"older continuous on $\Gamma$ and can be
calculated by Plemelj's formulas
\begin{equation}\label{e:plem}
\Phi^\pm (\eta(s))=\pm\frac{1}{2} \frac{\gamma(s) +\i\mu(s) }{
A(s)} +\frac{1}{2\pi \i} \int_J \frac{\gamma(t) +\i\mu(t) }{ A(t)}
\frac{\dot\eta(t)dt} {\eta(t)-\eta(s)}.
\end{equation}
The integral in~(\ref{e:plem}) is a Cauchy principal value
integral. It follows from~(\ref{e:plem}) that the boundary  values
satisfy the jump relation
\begin{equation}\label{e:plem-j}
A(s)\Phi^+(\eta(s))-A(s)\Phi^-(\eta(s)) = \gamma(s) + \i\mu(s).
\end{equation}

\noindent {\bf Interior RH problem:} Search a function $f$
analytic in $G$ with $f(\infty)=0$, continuous on $\Gamma\cup G$,
such that  the boundary values $f^+$ satisfy on $\Gamma$
\begin{equation}
\label{e:rhp} \Re [Af^+]=\gamma.
\end{equation}

To extend the results of~\cite{wegg,wegm} to unbounded multiply
connected regions, we shall borrow some definitions and notations
from~\cite{wegg,wegm}.
We define the following boundary value problem in the exterior
region $G^-$ as the exterior RH problem.

\noindent {\bf Exterior RH problem:} Search a function $g$
analytic in $G^-$, continuous on $\Gamma\cup G^-$, such that  the
boundary values $g^-$ satisfy on $\Gamma$
\begin{equation}
\label{e:erhp} \Re [Ag^-]=\gamma .
\end{equation}

We define the range spaces $R^\pm$ as the spaces of all real
functions $\gamma\in H$ for which the RH problems are solvable
and  the spaces $S^\pm$ as the spaces of the boundary values of solutions of the homogeneous RH problems, i.e.,
\begin{eqnarray}
  \label{e:R+}
  R^+ &:=& \{ \gamma \in H:\; \gamma =\Re [Af^+],\; f \mbox{
  analytic in $G$},\; f(\infty)=0\},  \\
  \label{e:S+}
  S^+ &:=& \{ \gamma \in H :\;\gamma=Af^+,\;f \mbox{ analytic in
  $G$},\;f(\infty)=0\}, \\
  \label{e:R-}
  R^- &:=& \{ \gamma \in H:\, \gamma=\Re [Ag^-],\; g \mbox{
analytic in $G^-$}\},\\
  \label{e:S-}
  S^- &:=& \{ \gamma \in H:\gamma=Ag^-,\, g \mbox{ analytic
in $G^-$}\}.
\end{eqnarray}
%

%--------------------------------------------------------------
\section{The generalized Neumann kernel}
\label{sc:gnk}
%--------------------------------------------------------------

%
We define the real  kernels $M$ and $N$ as real and imaginary parts (see~\cite{wegg,wegm} for details)
\begin{equation}
\label{e:M-N} M(s,t)+\i N(s,t):= \frac{1}{\pi}
\frac{A(s)}{A(t)}\frac{\dot\eta(t)}{\eta(t)-\eta(s)}.
\end{equation}
The kernel $N(s,t)$ is called the {\it generalized Neumann kernel}
formed with $A$ and~$\eta$.
It is continuous with
\begin{equation}
\label{e:N-d} N(t,t)= \frac{1}{\pi} \Im\left(\frac{1}{2}
\frac{\ddot\eta(t)}{\dot \eta(t)} - \frac{\dot A(t)}{ A(t)}
\right).
\end{equation}
When  $s,t\in J_k$  are in the same parameter interval $J_k$ then
\begin{equation}
\label{e:M-M1} M(s,t)= -\frac{1}{2\pi} \cot \frac{s-t}{2} +
M_1(s,t)
\end{equation}
with a continuous kernel $M_1$ which takes on the diagonal the
values
\begin{equation}\label{e:tM-d}
M_1(t,t)= \frac{1}{\pi}\Re \left(\frac{1}{2}
\frac{\ddot\eta(t)}{\dot \eta(t)} - \frac{\dot A(t)}{ A(t)}
\right).
\end{equation}
We define integral operators with the kernels $N$ and $M$ by
\begin{eqnarray}
  \label{e:opN}
 \bN \mu(s) &:=& \int_J N(s,t) \mu(t) dt, \quad s\in J, \\
 \label{e:opM}
 \bM \mu(s) &:=& \int_J M(s,t) \mu(t) dt, \quad s\in J.
\end{eqnarray}
The operator $\bN$ is a Fredholm integral operator and $\bM$ is a
singular operator where the integral in~(\ref{e:opM}) is a
principal value integral.
\begin{lemma}
\label{l:bdPhi} The boundary values of the function $\Phi$ defined
in~{\rm(\ref{e:Phi})} can be represented in terms of the operators
$\bN$ and $\bM$ by
\begin{equation}
 \label{e:Phi-bd-r}
2A\Phi^\pm = (\pm\bI+\bN-\i\bM)(\gamma+\i\mu).
\end{equation}
\end{lemma}
\noindent{\bf Proof:}
By multiplying both sides of~(\ref{e:plem}) by $2A(\zeta)$ and
using the definitions of the operators $\bN$ and $\bM$, we
obtain~(\ref{e:Phi-bd-r}).
\hfill $\Box$ \\
\begin{lemma}\label{l:bdv-anl}
{\rm (a)} If $f$ is an analytic function in the unbounded region $G$
with $\mbox{$f(\infty)=0$}$, then
\begin{equation}\label{e:bdv+anl}
(\bI-\bN+\i\bM)(Af^+) = 0.
\end{equation}
{\rm (b)} If $g$ is an analytic function in the bounded region $G^-$,
then
\begin{equation}\label{e:bdv-anl}
(\bI+\bN-\i\bM)(Ag^-) = 0.
\end{equation}
\end{lemma}
\noindent{\bf Proof:}%
(a) Let $\gamma := \Re[Af^+]$, $\mu := \Im[Af^+]$ and $\Phi$ be
formed with $\gamma,\,\mu$  according to~(\ref{e:Phi}).
By the Cauchy integral formula, we have $\Phi=f$ in $G$, which
implies that,
$A\phi^+= Af^+= \gamma+\i\mu$. Hence,~(\ref{e:bdv+anl}) follows
from~(\ref{e:Phi-bd-r}).

(b) Let $\gamma := \Re[Ag^-]$, $\mu := \Im[Ag^-]$ and $\Phi$ be
formed with $\gamma,\,\mu$  according to~(\ref{e:Phi}).
Since $G^-$ is on the right of $\Gamma$, the Cauchy integral
formula implies that $\Phi=-g$ in $G^-$.
Hence, $A\phi^-=-Ag^-= -(\gamma+\i\mu)$ and~(\ref{e:bdv-anl})
follows from~(\ref{e:Phi-bd-r}).
\hfill $\Box$ \\

There is a close connection between  RH problems and integral
equations with the generalized Neumann kernel (see
also~\cite{nasr,wegg,wegm}).

\begin{theorem}\label{t:rhp-ie}
If $f$ is a solution of the RH problem~{\rm(\ref{e:rhp})} with
boundary values
\begin{equation}\label{e:f-bd}
Af^+=\gamma+\i\mu,
\end{equation}
then the imaginary part $\mu$ in~{\rm(\ref{e:f-bd})} satisfies the
integral equation
\begin{equation}\label{e:ie}
\mu - \bN \mu =-\bM \gamma.
\end{equation}
\end{theorem}
\noindent{\bf Proof:}
Substituting~(\ref{e:f-bd}) into~(\ref{e:bdv+anl}) then taking the
imaginary part, we obtain~(\ref{e:ie}).
\hfill $\Box$ \\

It follows from the previous theorem that a solution of the RH
problem~(\ref{e:rhp}) yields a solution of the integral
equation~(\ref{e:ie}).
To use the integral equation~(\ref{e:ie}) to solve the RH
problem~(\ref{e:rhp}), we have the following theorem.
\begin{theorem}\label{t:ie-rhp}
Let the real function $\gamma\in H$ be given, $\mu$ be a solution
of~{\rm(\ref{e:ie})} and $\Phi$ be formed with $\gamma,\mu$
by~{\rm(\ref{e:Phi})}. Then $f:=\Phi$ in  $G$ satisfies
\begin{equation}
\label{e:bd-ft} Af^+=\gamma+h+\i\mu
\end{equation}
with
\begin{equation}\label{e:h}
    h = [\bM\mu-(\bI-\bN)\gamma]/2 \in  S^-.
\end{equation}
\end{theorem}
\noindent{\bf Proof:}
The jump relation~(\ref{e:plem-j}) implies that the function
$f:=\Phi$ in $G$ has the boundary value
\[
Af^+ = A\Phi^+ = A\Phi^- +\gamma +\i\mu.
\]
Since $\mu$ satisfies the integral equation~(\ref{e:ie}), it
follows from~(\ref{e:Phi-bd-r}) that
\[
2 A\Phi^- = -\gamma+\bN\gamma +\bM\mu
\]
which implies that $h:=A\Phi^-$ is real and given by~(\ref{e:h}).
Hence, $h \in S^-$.
\hfill $\Box$
%

%------------------------------------------------------------------
\section{The solvability}
\label{sc:slv}
%------------------------------------------------------------------

The solvability of RH problems on bounded multiply connected
regions was discussed in the classical text book~\cite{vek} (see
also~\cite{gak,wen}).
The solvability on unbounded multiply connected regions can be
deduced from the solvability on bounded regions by means of
M\"obius transform (see e.g.,~\cite[p.~141]{weg01}).
Let $z_0$ be a fixed point in $G^-$, say $z_0\in G_m$. The
unbounded region $G$ is transformed by means of the M\"obius 
transform
\begin{equation}\label{e:mob}
\Psi(z) := \frac{1}{z-z_0},
\end{equation}
onto a bounded multiply connected region $\hat G:= \Psi(G)$ of
connectivity $m$.
The M\"obius transform $\Psi$ transforms also the bounded exterior
region $G^-$ onto an unbounded region $\hat G^-$ exterior to the
boundary $\hat\Gamma := \partial \hat G=\Phi(\Gamma)$.
The boundary $\hat\Gamma$ is given by $\hat\Gamma := \hat\Gamma_0
\cup \hat\Gamma_1 \cdots\cup \hat\Gamma_{m-1}$ where $\hat\Gamma_0
:= \Psi(\Gamma_m)$ is the outer curve and is a counterclockwise
oriented; and the other curves $\hat\Gamma_j := \Psi(\Gamma_j)$,
$j=1,2,\ldots,m-1$, are clockwise oriented and are inside
$\hat\Gamma_0$.
The curve $\hat\Gamma$ is parametrized by
\begin{equation}\label{e:zeta-s}
\zeta(s) = \frac{1}{\eta(s)-z_0}, \quad s\in J.
\end{equation}
Let the function $\hat A$ be defined by
\begin{equation}\label{e:hA}
\hat A(s):=\zeta(s) A(s), \quad s\in J.
\end{equation}
Then, for a given function $\gamma\in H$, we define the RH
problems in the new regions $\hat G$ and $\hat G^-$ as follows.
\noindent {\bf Interior RH problem:} Search a function $\hat f$
analytic in $\hat G$, continuous on $\hat G\cup\hat\Gamma$, such
that the boundary values $\hat f^+$ satisfy on $\hat\Gamma$
\begin{equation}\label{e:rhp-h}
\Re [\hat A\hat f^+]=\gamma.
\end{equation}

\noindent {\bf Exterior RH problem:} Search a function $\hat g$
analytic in $\hat G^-$  with $\hat g(\infty)=0$, continuous on
$\hat G^-\cup\hat\Gamma$, such that  the boundary values $\hat
g^-$ satisfy on $\hat\Gamma$
\begin{equation}\label{e:erhp-h}
\Re [\hat A\hat g^-]=\gamma .
\end{equation}

We define also the spaces
\begin{eqnarray}
  \label{e:R+h}
  \hat R^+ &:=& \{ \gamma \in H:\; \gamma =\Re [\hat A\hat f^+],\; \hat f \mbox{
  analytic in $\hat G$}\}, \\
  \label{e:S+h}
  \hat S^+ &:=& \{ \gamma \in H :\;\gamma=\hat A\hat f^+,\;\hat f \mbox{ analytic in
  $\hat G$}\},\\
  \label{e:R-h}
  \hat R^- &:=& \{ \gamma \in H:\, \gamma=\Re [\hat A\hat g^-],\; \hat g \mbox{
  analytic in $\hat G^-$},\;\hat g(\infty)=0\}, \\
  \label{e:S-h}
  \hat S^- &:=& \{ \gamma \in H:\;\gamma=\hat A\hat g^-,\, \hat g \mbox{ analytic
  in $\hat G^-$},\; \hat g(\infty)=0\}.
\end{eqnarray}
\begin{lemma}\label{L:rhp-ub-b}
{\rm(a)} A function $f$ is a solution of the RH
problem~$(\ref{e:rhp})$ in the unbounded region $G$ if and only if
the function
\begin{equation}\label{e:fo-f}
\hat f(w) := \frac{(f\circ\Psi^{-1})(w)}{w}
\end{equation}
is a solution of the RH problem~$(\ref{e:rhp-h})$ in the bounded
region $\hat G$. \\
{\rm(b)} A function $g$ is a solution of the exterior RH
problem~$(\ref{e:erhp})$ in the bounded region $G^-$ if and only
if the function
\begin{equation}\label{e:go-g}
\hat g(w) := \frac{(g\circ\Psi^{-1})(w)}{w}
\end{equation}
is a solution of the RH problem~$(\ref{e:erhp-h})$ in the
unbounded region $\hat G^-$.
\end{lemma}
\noindent{\bf Proof:}
The proof follows from the definitions of the M\"obius transform
$\Psi$, the function $\hat A$ and the RH
problems~(\ref{e:rhp}),~(\ref{e:erhp}),~(\ref{e:rhp-h}),~(\ref{e:erhp-h}).
\hfill $\Box$ \\
\begin{lemma}\label{L:sp-sp-h}
The spaces $R^{\pm}$, $S^{\pm}$, $\hat R^{\pm}$ and $\hat S^{\pm}$
satisfy
\begin{equation}\label{e:sp-sp-h}
\hat R^+ = R^+, \quad \hat S^+ = S^+, \quad \hat R^- = R^-, \quad
\hat S^- = S^-.
\end{equation}
\end{lemma}
\noindent{\bf Proof:}
In view of Lemma~\ref{L:rhp-ub-b}, the proof follows from the
definitions of the spaces $R^{\pm}$, $S^{\pm}$, $\hat R^{\pm}$ and
$\hat S^{\pm}$.
\hfill $\Box$ \\
%

%%
%The solvability of RH problems as well as integral equations with
%the generalized Neumann kernel depends on the index of the
%function $A$ on the boundary $\Gamma$.
%%
The index $\kappa_j$ of the function $A$ on the curve $\Gamma_j$
is defined as the winding number of  $A$ with respect to $0$,
\begin{equation}\label{e:ind-j}
\kappa_j := \frac{1}{2\pi}\Delta\arg(A)|_{\Gamma_j}, \quad j=1,2,
\ldots,m,
\end{equation}
i.e., the change of the argument of $A$ along the curve $\Gamma_j$
divided by $2\pi$.  The index $\kappa$ of the function $A$ on the
whole boundary curve $\Gamma$ is the sum
\begin{equation}\label{e:ind}
\kappa := \sum_{j=1}^m \kappa_j.
\end{equation}
The index $\hat\kappa_{j}$ of the function $\hat A$ on the curve
$\hat\Gamma_j$ and the index $\hat\kappa$ of the function $\hat A$
on the whole boundary $\hat\Gamma$ are easily calculated from the
index of the function $A$ by
\begin{equation}\label{e:ind-hA}
\hat\kappa_{0} = \kappa_m+1, \quad \hat\kappa_{j}=\kappa_j, \quad
j=1,2,\ldots,m-1, \quad \hat\kappa =\kappa+1.
\end{equation}
The RH problems~(\ref{e:rhp-h}) and~(\ref{e:erhp-h}) are of the
types studied in~\cite{wegm}.
Thus, in view of Lemma~\ref{L:sp-sp-h} and Eq.~(\ref{e:ind-hA}),
we have the following results from~\cite{wegm} for the solvability
of the RH problems~(\ref{e:rhp}) and~(\ref{e:erhp}).
\begin{lemma}[\cite{wegm}]
The spaces $R^-$ and $S^\pm$ satisfy
\begin{eqnarray}
  \label{e:S-S+}
  S^-\cap S^+ &=& \{0\}, \\
  \label{e:R-S+}
  R^-\cap S^+ &=& \{0\}.
\end{eqnarray}
\end{lemma}

\begin{theorem}[\cite{wegm}]\label{t:dim-S}
The dimension of the space $S^-$ and the codimension of the space
$R^-$ are determined by the index of $A$ as follows:
\begin{eqnarray}
  \label{e:dim-S}
\dim(S^-) &=& \sum_{j=1}^m\max(0,2\kappa_j+1), \\
  \label{e:cdim-R}
\codim(R^-) &=& \sum_{j=1}^m\max(0,-2\kappa_j-1).
\end{eqnarray}
\end{theorem}

\begin{theorem}[\cite{wegm}]\label{t:dim+S}
The dimension of the space $S^+$ and the codimension of the space
$R^+$ are determined by the index of $A$ as follows:\\
{\rm(a)} If $\kappa\ge0$, then
\begin{equation}\label{e:dim+S1}
    \dim(S^+)=0, \quad \codim(R^+) = 2\kappa+m.
\end{equation}
{\rm(b)} If $\kappa\le -m$, then
\begin{equation}\label{e:dim+S2}
    \dim(S^+)=-2\kappa-m, \quad \codim(R^+) = 0.
\end{equation}
{\rm(c)} If $-m+1\le\kappa\le-1$, then
\begin{equation}\label{e:dim+S3}
    -2\kappa-m \le\dim S^+\le -\kappa, \quad 2\kappa+m\le\codim(R^+) \le m+\kappa.
\end{equation}
\end{theorem}
%
%\noindent{\bf Proof:}
%%
%In view of Lemma~\ref{L:sp-sp-h} and Eq.~(\ref{e:ind-hA}), the
%proof follows from~\cite[Theorem~9]{wegm}.
%%
%\hfill $\Box$ \\
%

%
For studying the solvability of the integral
equation~(\ref{e:ie}), it follows from~(\ref{e:zeta-s}) that
\[
\eta(s) = \frac{1}{\zeta(s)} + z_0, \quad \dot\eta(s) =
-\frac{\dot\zeta(s)}{\zeta(s)^2}, \quad s\in J.
\]
Hence
\begin{equation}\label{e:A-tA}
\frac{A(s)}{A(t)} \frac{\dot\eta(t) }{\eta(t) - \eta(s)} = \frac{\hat A(s)}{\hat
A(t)} \frac{\dot\zeta(t)}{\zeta(t)-\zeta(s)}
\end{equation}
where $\hat A$ is defined by~(\ref{e:hA}). Let the real kernels
$\hat M$ and $\hat N$ be defined by
\begin{equation}
 \label{e:hN-hM}
\hat M(s,t) +\i\hat N(s,t) := \frac{1}{\pi}  \frac{\hat A(s)}{\hat
A(t)}\frac{\dot\zeta(t)}{\zeta(t)-\zeta(s)}.
\end{equation}
i.e., the kernel $\hat N$ is the generalized Neumann kernel formed
with $\hat A$ and $\zeta$.
Let $\hat\bM$ and $\hat\bN$ be the integral operators defined on
$H$ with the kernels $\hat M$ and $\hat N$.
In view of~(\ref{e:A-tA}), we have $\hat N(s,t) = N(s,t)$ and $\hat M(s,t) = M(s,t)$ for all $(s,t) \in J\times J$.
Hence
\begin{equation}\label{e:op-h}
\hat \bN = \bN, \quad \hat \bM = \bM.
\end{equation}
The operators $\hat\bM$ and $\hat\bN$ are of the types studied
in~\cite{wegm}.
Thus, in view of Lemma~\ref{L:sp-sp-h} and Eq.~(\ref{e:ind-hA}),
we have the following results from~\cite{wegm} for the properties
of the integral operators $\bN$ and $\bM$.
\begin{lemma}[\cite{wegm}]
The operators $\bN,\bM$ and the identity operator $\bI$ are
connected by the following relations:
\begin{eqnarray}
  \bN^2-\bM^2 &=& \bI, \\
  \bN\bM+\bM\bN &=& 0.
\end{eqnarray}
\end{lemma}
%
%\noindent{\bf Proof:}
%%
%The lemma is proved with the same approach used
%in~\cite[Lemma~4]{wegm}.
%%
%\hfill $\Box$ \\
%

%
\begin{theorem}[\cite{wegm}]
The range-spaces and the null-spaces of the operators $\bM$ and
$\bI\pm \bN$ are related to the spaces $S^\pm$ and $R^\pm$ by
\begin{eqnarray}
  \label{e:rang-M}
\Range(\bM) &=& R^+\cap R^-,  \\
  \label{e:null-M}
\Null(\bM) &=& S^+\oplus S^-,  \\
 \label{e:rang+N}
\Range(\bI+\bN) &\subset& R^+,  \\
 \label{e:null+N}
\Null(\bI+\bN) &=& S^-,  \\
 \label{e:rang-N}
\Range(\bI-\bN)&=&R^-, \\
 \label{e:null-N}
\Null(\bI-\bN) &=& S^+\oplus W,
\end{eqnarray}
where $W$ is isomorphic $($via $\bM)$ to $R^+\cap S^-$.
\end{theorem}
%
%\noindent{\bf Proof:}
%%
%In view of Lemma~\ref{L:sp-sp-h} and Eq.~(\ref{e:op-h}):\\
%The Formula~(\ref{e:rang-M}) follows
%from~\cite[Formula~(71)]{wegm}. \\
%The Formula~(\ref{e:null-M}) follows
%from~\cite[Formula~(56)]{wegm}. \\
%The Formula~(\ref{e:rang+N}) follows
%from~\cite[Formula~(69)]{wegm}. \\
%The Formula~(\ref{e:null+N}) follows
%from~\cite[Formula~(102)]{wegm}. \\
%The Formula~(\ref{e:rang-N}) follows
%from~\cite[Formula~(103)]{wegm}. \\
%The Formula~(\ref{e:null-N}) follows
%from~\cite[Lemma~20(a)]{wegm}.
%%
%\hfill $\Box$ \\
%

%
\begin{corollary}\label{C:sl-ie}
The integral equation~{\rm(\ref{e:ie})} is solvable for all
$\gamma\in H$.
\end{corollary}
\noindent{\bf Proof:}
The Formulas~(\ref{e:rang-M}) and~(\ref{e:rang-N}) imply that
\[
\Range(\bM) \subset \Range(\bI-\bN),
\]
which implies that the integral equation~{\rm(\ref{e:ie})} is
always solvable.
\hfill $\Box$ \\
%

%%
%The exact number of linearly independent solutions of the
%homogeneous integral equation is determined by the index as in the
%following theorem.
%
\begin{theorem}[\cite{wegm}]\label{t:dim-+N}
The dimensions of the null-spaces of the operators $\bI\pm\bN$ are
given by
\begin{eqnarray}
 \label{e:dim+N}
\dim(\Null(\bI+\bN)) &=&\sum_{j=1}^m\max(0,2\kappa_j+1),  \\
 \label{e:dim-N}
\dim(\Null(\bI-\bN)) &=& \sum_{j=1}^m\max(0,-2\kappa_j-1).
\end{eqnarray}
\end{theorem}
%
%\noindent{\bf Proof:}
%%
%In view of Eqs.~(\ref{e:ind-hA}) and~(\ref{e:op-h}):\\
%The Formula~(\ref{e:dim+N}) follows
%from~\cite[Formula~(105)]{wegm}. \\
%The Formula~(\ref{e:dim-N}) follows
%from~\cite[Formula~(106)]{wegm}.
%%
%\hfill $\Box$ \\
%

%
\begin{lemma}[\cite{wegm}]
If $\lambda$ is an eigenvalue of $\bN$ with eigenfunction
$\nu\notin S^+\cup S^-$, then $-\lambda$ is also an eigenvalue of
$\bN$ and $\bM\nu$ is a corresponding eigenfunction.
\end{lemma}
%
%\noindent{\bf Proof:}
%%
%The lemma is proved with the same approach used
%in~\cite[Lemma~10]{wegm}.
%%
%\hfill $\Box$ \\
%

It follows from Theorem~\ref{t:dim+S} that the RH
problem~(\ref{e:rhp}) is not necessary solvable for general $A$
and $\gamma$.
In the following theorem, we shall show that the right-hand side
of~(\ref{e:rhp}) can be always modified such that the new problem
is solvable.
\begin{theorem}\label{T:rhp-p}
For any $\gamma\in H$, there exists a real function $h\in S^-$
such that the RH problem
\begin{equation}\label{e:rhp-p}
    \Re [Af^+] = \gamma+h
\end{equation}
is solvable.
\end{theorem}
\noindent{\bf Proof:}
For any $\gamma\in H$, it follows from Corollary~\ref{C:sl-ie}
that the integral equation~(\ref{e:ie}) is solvable.
Let $\mu$ be a solution of the integral equation~(\ref{e:ie}).
Then, Theorem~\ref{t:ie-rhp} implies that a real function $h:=
[\bM\mu-(\bI-\bN)\gamma]/2\in S^-$ exists such that $Af^+ =
\gamma+h+\i\mu$ are boundary values of an analytic solution $f$ in
$G$. Hence, $f$ is a solution of the RH problem~(\ref{e:rhp-p}).
\hfill $\Box$ \\
\begin{theorem}\label{T:H=S+R}
If $\dim(S^+) = \dim(\Null(\bI-\bN))$, then the space $H$ has the
decomposition
\begin{equation}\label{e:H=S+R}
H=R^+\oplus S^-.
\end{equation}
\end{theorem}
\noindent{\bf Proof:} %
Let $\gamma\in H$ be a given function. It follows from
Theorem~\ref{T:rhp-p} that a function $h\in S^-$ exists such that
$\gamma+h\in R^+$. Hence
\begin{equation}\label{e:dim-SW}
H= R^+ + S^-.
\end{equation}
Since $\dim(S^+) = \dim(\Null(\bI- \bN))$, it follows
from~(\ref{e:null-N}) that
\[
\dim(R^+\cap S^-)=\dim(\Null(\bI-\bN))-\dim(S^+)=0,
\]
which implies that $R^+ \cap S^- =\{0\}$.
Hence, the sum in~(\ref{e:dim-SW}) is direct.
\hfill $\Box$ \\

\begin{corollary}\label{C:H=S+R}
Let $\dim(S^+) = \dim(\Null(\bI-\bN))$. Then for any $\gamma\in H$, there
exists a unique function $h\in S^-$ such that the RH
problem~$(\ref{e:rhp-p})$ is solvable.
\end{corollary}
%

%-------------------------------------------------------------
\section{The Modified Dirichlet problem}
\label{sc:mdir}
%------------------------------------------------------------------

In this section, we shall study the solvability of the RH
problem~(\ref{e:rhp}) and the integral equation~(\ref{e:ie}) for
the special case
\begin{equation}\label{e:A=1}
A=1.
\end{equation}
In this case, the kernel $N$ is known as the Neumann kernel~\cite[p.~286]{hen} and the RH problem~(\ref{e:rhp}) is known as the modified
Dirichlet problem~\cite{gak,mik,mus,wen} or the Schwartz
problem~\cite{gak}.
This special case is of practical use in conformal mapping of unbounded multiply connected regions (see~\cite{nasm}).

The function $A=1$ has the index
\begin{equation}\label{e:ind-A=1}
\kappa_j=0, \quad j=1,2,\ldots,m, \quad \kappa=0.
\end{equation}
Hence, Theorems~\ref{t:dim-S},~\ref{t:dim+S} and~\ref{t:dim-+N}
imply that
\begin{equation}\label{e:S+A=1}
\dim(S^-)=m, \qquad \dim(S^+) =0, \qquad \dim(\Null(\bI-\bN)) = 0.
\end{equation}
Let the real function $\chi^{[j]}(s)$, $j=1,2, \ldots,m$, be
defined for $s\in J$ by
\begin{equation}\label{e:chi}
\chi^{[j]}(s):= \left\{ %
\begin{array}{l@{\hspace{0.5cm}}l} %
1,  & s\in J_j,\\   %
0,  & s\notin J_j,  %
\end{array}
\right .
\end{equation}
and the sectionally analytic function $g_j(z)$ be defined by
\begin{equation}\label{e:g-k}
g_j(z):= \left\{ %
\begin{array}{l@{\hspace{0.5cm}}l} %
1,  & z\in G_j,\\   %
0,  & \mbox{$z$ in the exterior domain of the curve $\Gamma_j$}. %
\end{array}
\right .
\end{equation}
Hence $\chi^{[j]} = Ag_j^-\in S^-$ for all $j = 1,2,\ldots,m$.
From~(\ref{e:S+A=1}), we have $\dim(S^-)=m$.
Since $\chi^{[1]},\chi^{[2]}, \ldots, \chi^{[m]}$ are linearly
independent, we obtain
\begin{equation}\label{e:S-A=1}
S^- = {\rm span} \{\chi^{[1]}, \chi^{[2]}, \ldots, \chi^{[m]}\}.
\end{equation}
\begin{theorem}\label{T:mdir}
For any $\gamma\in H$, there exists a unique solution $\mu$ of the
integral equation~$(\ref{e:ie})$ and a unique function $h =
(h_1,h_2, \ldots,h_m)$ given by~$(\ref{e:h})$ such that
$f^+=\gamma+h+\i\mu$ are boundary values of the unique solution of
the RH problem~$(\ref{e:rhp-p})$.
\end{theorem}
\noindent {\bf Proof.} %
By~(\ref{e:S+A=1}), $\Null(\bI-\bN)=\{0\}$ which implies, in view
of Fredholm alternative theorem, that the integral
equation~(\ref{e:ie}) has a unique solution $\mu$.
Then, Theorem~\ref{t:ie-rhp} implies that $f^+=\gamma+h+\i\mu$ are
boundary values of an analytic function $f$ in $G$ where $h\in
S^-$ is given by~{\rm(\ref{e:h})}.
Since $\dim(S^+)=0$, $f$ is the unique solution of the RH
problem~(\ref{e:rhp-p}).
In view of Eq.~(\ref{e:S+A=1}), Corollary~\ref{C:H=S+R} and
Eq.~(\ref{e:S-A=1}) imply that the function $h$ is unique and $h =
(h_1,h_2, \ldots,h_m)$ with real constants $h_1,h_2, \ldots,h_m$.
\hfill $\Box$
%

%-----------------------------------------------
\section{Conclusions}
\label{sc:con}

%This paper extended the results obtained in the previous
%papers~\cite{wegg,wegm} to unbounded multiply connected regions.
%
We have derived and studied a boundary integral equation with the
generalized Neumann kernel for the RH problem on unbounded
multiply connected regions.
By means of a M\"obius transform, the solvability of the derived
integral equation was obtained from the related known results for
bounded regions.
The derived integral equation was then used to obtain boundary
integral equations for the modified Dirichlet problem on
unbounded multiply connected regions.
The boundary integral equations were derived in this paper for
regions with smooth boundaries.
Nevertheless, these equations, with slight modifications, can be
applied to regions with corners (see~\cite{nasr} for
simply connected regions case).
Several accurate numerical methods are available for solving
boundary integral equations (see e.g.~\cite{atk}).
For regions with smooth boundary, one can use the Nystr\"om method
with the trapezoidal rule~\cite{atk}.
For regions with corners, we can use a Nystr\"om method based on
the trapezoidal rule with a graded mesh~\cite{atk,kre}.
%
%
%No numerical results were presented in this paper.
%
%Numerical solution of RH problems on multiply connected regions by
%integral equations with the generalized Neumann kernel will be
%given in future papers.
%

\end{document}